\documentclass[twocolumn]{autart}    

\usepackage{amsmath,amssymb,amsfonts,bm}
\usepackage{graphicx}
\usepackage{algorithm}
\usepackage{algorithmicx,algpseudocode}
\usepackage{cite}
\usepackage{color}
\let\theoremstyle\relax
\usepackage{amsthm}
\usepackage{cuted}
\usepackage[mathscr]{euscript}
\theoremstyle{plain}
\newtheorem{theorem}{Theorem}[section]

\newtheorem{proposition}[theorem]{Proposition}
\theoremstyle{definition}

\theoremstyle{remark}

\usepackage{wrapfig}

\begin{document}

\begin{frontmatter}

\title{General Discrete-Time Fokker-Planck Control by Power Moments} 

\author[Guangyu]{Guangyu Wu}\ead{chinarustin@sjtu.edu.cn},  
\author[Anders]{Anders Lindquist}\ead{alq@kth.se}              

\address[Guangyu]{Department of Automation, Shanghai Jiao Tong University, Shanghai, China}  
\address[Anders]{Department of Automation and School of Mathematical Sciences, Shanghai Jiao Tong University, Shanghai, China}
          

\begin{abstract}
In this paper, we address the so-called general Fokker-Planck control problem for discrete-time first-order linear systems. Unlike conventional treatments, we don't assume the distributions of the system states to be Gaussian. Instead, we only assume the existence and finiteness of the first several order power moments of the distributions. It is proved in the literature that there doesn't exist a solution, which has a form of conventional feedback control, to this problem. We propose a moment representation of the system to turn the original problem into a finite-dimensional one. Then a novel feedback control term, which is a mixture of a feedback term and a Markovian transition kernel term is proposed to serve as the control input of the moment system. The states of the moment system are obtained by maximizing the smoothness of the state transition. The power moments of the transition kernels are obtained by a convex optimization problem, of which the solution is proved to exist and be unique. Then they are mapped back to the probability distributions. The control inputs to the original system are then obtained by sampling from the realized distributions. Simulation results are provided to validate our algorithm in treating the general discrete-time Fokker-Planck control problem.
\end{abstract}

\end{frontmatter}

\section{Introduction}
\label{sec:introduction}

Distribution steering is an important problem in controls and robotics, of which the significance is underscored by its relevance to various real-world scenarios, including quality control of manufacturing process, formation control of swarm robots, etc. Typically, the previous works on distribution steering can be separated into two types of problems, which we call as the Liouville control problem and the Fokker-Planck control problem. The Liouville control problem considers the system equation to be linear differential/difference equations without noise. For the Fokker-Planck control problem, an additive noise, which is a Wiener process for the continuous-time case or a zero-mean Gaussian random variable for the discrete-time case, is taken into consideration. 

We first review the relevant results on the two types of distribution steering problems. The continuous-time Liouville control problem, involving perfect state information, was originally introduced by Brockett \cite{brockett2012control}. In his work, Brockett \cite{brockett2012control} delves into the task of steering the state of a continuous-time linear system. This system's initial state follows a known (multivariate) Gaussian distribution, and the objective is to drive it towards a desired (multivariate) Gaussian distribution by a predetermined finite terminal time. Moreover, conditions that determine the feasibility of solving this steering problem are established. Considering the case of static output feedback and parallel interconnections of linear systems, an extension of the results in \cite{brockett2012control} is given in \cite{dirr2016controlling}. A more recent treatment to the Liouville control problem of the discrete-time SISO linear systems is proposed in \cite{bakolas2019dynamic}. In the literature, there are several topics of research, which have similar tasks as Liouville control, however with different names. Examples include ensemble control \cite{li2007ensemble, qi2013optimal, li2010ensemble} and control of families of systems \cite{fuhrmann2015mathematics}. 

In the Liouville control problem, it is assumed that there exists no noise in the system equation. However, the landscape of real-world engineering practices introduces an inescapable reality: the presence of both external and internal noises, encompassing elements such as modeling errors. This naturally propels us to delve into the consideration of noise within the system equation itself. In the realm of the Fokker-Planck equation, a notable facet emerges – the inclusion of an additive Gaussian noise term. Pertinent investigations into the Fokker-Planck control problem are well-documented in existing literature. Expanding the horizon to encompass continuous-time Fokker-Planck equations, which fundamentally correspond to Schrödinger bridge problems, prior endeavors have been documented in works such as \cite{blaquiere1992controllability, pra1990markov, pavon1991free, filliger2008connection, chen2016optimal, chen2015optimal, chen2015optimal2, chen2018optimal}. Notably, discrete-time Fokker-Planck control problem explorations are also prevalent, as exemplified by \cite{goldshtein2017finite, okamoto2018optimal, balci2020covariance, balci2021covariance, balci2023covariance, pilipovsky2023data}, often referred to as covariance steering within the literature.

However, in the literature, the distributions of the system states are always assumed to be Gaussian, which hinders the use of the previous results in a wider range of engineering applications. Recent treatments have been proposed considering general types of noises, rather than constraining them to be Gaussian. In \cite{sivaramakrishnan2022distribution}, characteristic functions are utilized to treat the distribution steering under general disturbances. A martingale noise model is considered in \cite{liu2022optimal}, where the noise can be a Wiener process, a Wiener process with jumps in it, or a mixture of them. This treatment broadens the conventional assumption on the system noise, which makes it more applicable in real engineering problems. In our previous works \cite{wu2022density, wu2022group, wu2023general}, we considered the general discrete-time Liouville control problem, where the distributions of the system states are not confined to be Gaussian. With very loose constraint on the distributions, our works can be well applied to the control of a very large arbitrary group of agents of which the dynamic equation is a Liouville one, which was the first time in the literature.  

In this paper, we propose to treat the general discrete-time Fokker-Planck control problem. This
type of control problem aims at steering an initial distribution to a terminal one, where the system dynamic is governed by a discrete-time Fokker-Planck equation. It is called "general" since there is no assumption on the initial and terminal distributions being falling into specific function classes (a conventional choice is the Gaussian distribution). Instead, we only assume the distributions to have first several orders of power moments. It is not a strict constraint, which is satisfied by most commonly seen distributions, except for distributions of which all the power moments don't exist, e.g. Cauchy distribution. By this problem setting, the general Fokker-Planck control problem becomes an infinite-dimensional problem. 

The structure of the paper unfolds as follows. In Section 2, we offer a comprehensive formulation of the overarching Fokker-Planck control problem we are set to tackle. Moving on to Section 3, our focus shifts to the pivotal task of reducing the dimension of the original infinite-dimensional problem. Here, we introduce a moment system representation tailored to the discrete-time Fokker-Planck equation. Within this framework, we propose a dynamic control law, which is a mixture of a feedback component and a Markovian transition kernel. A condition on the states of the moment system is proposed to ensure the existence of control inputs, given the discrete-time Fokker-Planck equation and the additive system noise. In Section 4, We propose a control scheme by convex optimization, which ensures the existence and uniqueness of the control inputs for the moment system at each time step, provided with the states of the moment system. Then a solution to the inverse problem of a map from the power moments to the infinite-dimensional transition kernel is proposed. An algorithms considering the general discrete-time Fokker-Planck control problem is proposed. Two numerical examples are given in Section 5, which validate the proposed algorithm in treating the discrete-time Fokker-Planck control problem. A conclusion of the paper and future research directions are given in Section 6.

\section{Problem formulation}
\label{sec:Problem}

We first introduce several notations that will be used throughout this paper. We denote the set of nonnegative integers and strictly positive integers as $\mathbb{N}$ and $\mathbb{N}_{+}$ respectively. The set of nonnegative integers from $a$ to $b$ is defined as
$$
\{a, \cdots, b\}=\left[ a, b \right] \cap \mathbb{N}:=\mathbb{N}_{a}^{b}.
$$
$\mathbb{E}\left[ \cdot \right]$ denotes the expectation operator. Moreover, we denote the set of all feasible probability distributions defined on $\mathbb{R}$, of which the first $2n$ orders of power moments exist and are finite, as $\mathcal{P}_{2n}$, namely
$$
\begin{aligned}
\mathcal{P}_{2n} := & \left\{ p(x)>0 \mid \int_{\mathbb{R}} p(x) dx = 1, \right. \\
& \left. \left|\int_{\mathbb{R}} x^{i}p(x) dx \right| < \infty, i \in \mathbb{N}_{0}^{2n} \right\}.
\end{aligned}
$$

Most commonly seen probability distributions fall within the set $\mathcal{P}_{2n}$, except for distributions of which the power moments don't exist, e.g. the Cauchy distribution. 

For a given $K \in \mathbb{N}_{+}$, let $\left\{a(k) \in \mathbb{R}: k \in \mathbb{N}_{0}^{K-1}\right\},\{b(k) \in$ $\left.\mathbb{R}: k \in \mathbb{N}_{0}^{K-1}\right\}$ denote known system parameters. Denote the system state, the control input and the additive white noise as $x(k)$, $u(k)$ and $w(k)$ respectively. Assume that $x(k)$, $u(k)$ and $w(k)$ are random variables, and their probability distributions are denoted as $p_{k}(x), q_{k}(u)$ and $h_{k}(w)$. Consider the following discrete-time first-order linear system
\begin{equation}
x(k+1) =a(k) x(k)+b(k) u(k)+w(k),
\label{syseq}
\end{equation}
for $k \in \mathbb{N}_{0}^{K-1}$. The Fokker-Planck control problem can then be formulated as follows. Given the system equation \eqref{syseq}, a sequence of independent Gaussian noises $\left\{ w(k): k \in \mathbb{N}_{0}^{K-1}\right\}$ and an initial state $x(0)$, which is a sample drawn from the prescribed initial distribution $p_{0}(x) \in \mathcal{P}_{2n}$, find a sequence of feasible control inputs $\left\{u(k): k \in \mathbb{N}_{0}^{K-1}\right\}$, which steers the terminal system state $x(K)$ to follow a desired distribution $\tau(x)$. In the conventional treatments of the distribution steering, the initial and terminal distributions are assumed to be Gaussian, and the control problem is reduced to steering the first and second order power moments to desired ones, by solving the algebraic Riccati or Sylvester equations. However, in our problem setting, we only assume the initial and terminal distributions to fall within the set $\mathcal{P}_{2n}$. By this setting, the problem is an infinite-dimensional one, which makes the original problem intractable and open. In other words, a proper dimension reduction method needs to be proposed to treat this problem. 

\section{A moment representation of the original system}
In our previous papers \cite{wu2022density, wu2022group, wu2023general}, we proposed to use power moments to decrease the dimension of the general Liouville control problem. However, to extend the results of these works to the Fokker-Planck control problem is not a trivial problem. Introducing the additive noise may cause some of the system states not reachable. In the following parts of the paper, we will propose a moment representation of the original problem, then obtain a control scheme for the moment system, with proves to the reachability of the states and existence of the corresponding control inputs. 

We first give a moment representation of the original system \eqref{syseq}. It has been proven in \cite{elamvazhuthi2018optimal} that the distribution steering problem, where the initial and terminal distributions belong to different function classes, cannot be solved using conventional deterministic feedback laws. Therefore, we propose a new type of feedback control law. It is a mixture of a feedback portion of the system state and a random variable, serving as the transition kernel of Markovian process, which reads
\begin{equation}
u(k) = -c(k)a(k)x(k) + \mathscr{F}(k), \quad 0 \leq c(k) \leq 1.
\label{control}
\end{equation}
Here $\mathscr{F}(k)$ is a random variable, inspired by the result in \cite{biswal2021decentralized}, which is independent of $\left\{ x(t): t \in \mathbb{N}_{0}^{k} \right\}$. By this choice of control input, we can write the original system equation as
$$
\begin{aligned}
& \mathbb{E}\left[u^{\ell}(k)\right] \\
= & \sum_{i=0}^{\ell}(-c(k)a(k))^i \mathbb{E}\left[x^i(k)\right] \mathbb{E}\left[\mathscr{F}^{\ell-i}(k)\right] .
\end{aligned}
$$
Furthermore, by instituting $u(k)$ into \eqref{syseq}, we obtain
\begin{equation}
\begin{aligned}
& x(k+1)\\
= & a(k)\left( 1 - b(k)c(k) \right)x(k) + b(k)\mathscr{F}(k) + w(k)\\
= & \tilde{a}(k)x(k) + \tilde{u}(k).
\label{newsyseq}
\end{aligned}
\end{equation}
where we have defined
$$
\tilde{a}(k) := a(k)\left( 1 - b(k)c(k) \right),
$$
and
\begin{equation}
\tilde{u}(k) := b(k)\mathscr{F}(k) + w(k).
\label{tildeu}
\end{equation}
Since $\mathscr{F}(k), w(k)$ are both independent of $x(k)$, we have that $\tilde{u}(k)$ is also independent of $x(k)$. The power moments of $x(k+1)$ can then be written as
$$
\mathbb{E}\left[ x^{\ell}(k+1) \right] = \sum_{i=0}^{\ell}\tilde{a}^{i}(k) \mathbb{E}\left[x^i(k)\right] \mathbb{E}\left[\tilde{u}^{\ell-i}(k)\right].
$$
Then the dynamics of the moments up to order $2n$ can be written as a linear matrix equation
\begin{equation}
\mathscr{X}(k+1)=\tilde{\mathscr{A}}(\tilde{\mathscr{U}}(k)) \mathscr{X}(k)+\tilde{\mathscr{U}}(k),
\label{momentsyseq}
\end{equation}
where the new system state reads
\begin{equation}
\mathscr{X}(k) = \begin{bmatrix}
\mathbb{E}[x(k)] & \mathbb{E}[x^{2}(k)] & \cdots & \mathbb{E}[x^{2n}(k)]
\end{bmatrix}^{\intercal},
\label{XK}
\end{equation}
and the new control input reads
\begin{equation}
\tilde{\mathscr{U}}(k) = \begin{bmatrix}
\mathbb{E}[\tilde{u}(k)] & \mathbb{E}[\tilde{u}^{2}(k)] & \cdots & \mathbb{E}[\tilde{u}^{2n}(k)]
\end{bmatrix}^{\intercal}.
\label{UK}
\end{equation}

The new system matrix $\tilde{\mathscr{A}}(\tilde{\mathscr{U}}(k))$ is given in \eqref{longeq1}. Note that $\tilde{\mathscr{A}}(\tilde{\mathscr{U}}(k))$ is different from the conventional system matrices, since it contains the control input $\tilde{u}(k)$. Then the task of determining the sequence of the control inputs $\left\{ u(k):k \in \mathbb{N}_{0}^{K} \right\}$ becomes one of determining $\left\{ c(k):k \in \mathbb{N}_{0}^{K} \right\}$ and $\left\{ \tilde{u}(k):k \in \mathbb{N}_{0}^{K} \right\}$. 

\begin{figure*}[t]
\begin{equation}
\tilde{\mathscr{A}}(\tilde{\mathscr{U}}(k))
= \begin{bmatrix}
\tilde{a}(k) & 0 & 0 & \cdots & 0\\ 
2\tilde{a}(k)\mathbb{E}[\tilde{u}(k)] & \tilde{a}^{2}(k) & 0 & \cdots & 0\\ 
3\tilde{a}(k)\mathbb{E}[\tilde{u}^{2}(k)] & 3\tilde{a}^{2}(k)\mathbb{E}[\tilde{u}(k)] & \tilde{a}^{3}(k) & \cdots & 0\\ 
\vdots & \vdots & \vdots & \ddots\\ 
\binom{2n}{1}\tilde{a}(k)\mathbb{E}[\tilde{u}^{2n-1}(k)] & \binom{2n}{2}\tilde{a}^{2}(k)\mathbb{E}[\tilde{u}^{2n-2}(k)] & \binom{2n}{3}\tilde{a}^{3}(k)\mathbb{E}[\tilde{u}^{2n-3}(k)] &  & \tilde{a}^{2n}(k)
\end{bmatrix}
\label{longeq1}
\end{equation}
\hrulefill
\vspace*{4pt}
\end{figure*}

However, we should note that there exists an implicit constraint on $\tilde{u}(k)$, which is different from the results in our works on the discrete Liouville control problems. By \eqref{tildeu}, we have that
\begin{equation}
\begin{aligned}
& \mathbb{E}\left[ \tilde{u}^{\ell}(k) \right]\\
= & \sum_{i=0}^{\ell}b^{i}(k)\mathbb{E}\left[\mathscr{F}^i(k)\right] \mathbb{E}\left[w^{\ell-i}(k)\right], \ell \in \mathbb{N}_{1}^{2n}
\label{Expecttildeul}
\end{aligned}
\end{equation}
where the equation is because $\mathscr{F}(k)$ and $w(k)$ are independent. Therefore, given the variance of the white noise $\sigma^{2}$, we can obtain its power moments up to order $2n$ by
$$
\mathbb{E}\left[w^{2\ell}(k)\right] = \sigma^{2\ell}(2\ell - 1)!!, \ell \in \mathbb{N}_{1}^{n}
$$
and
$$
\mathbb{E}\left[w^{2\ell-1}(k)\right] = 0, \ell \in \mathbb{N}_{1}^{n}.
$$
where $!!$ denotes the dual factorial. Then by \eqref{Expecttildeul}, we can calculate $\mathbb{E}\left[\mathscr{F}^i(k)\right]$ for $i \in \mathbb{N}_{1}^{2n}$. Define the Hankel matrix of a random variable $\eta(k)$ as $H_{\eta}(k)$. The condition of existence of $\mathscr{F}(k)$ is that the Hankel matrix of it is positive semidefinite, i.e., 
\begin{equation}
\begin{aligned}
& H_{\mathscr{F}}(k)\\
= & \begin{bmatrix}
1 & \mathbb{E}\left[ \mathscr{F}(k) \right ] & \cdots & \mathbb{E}\left[ \mathscr{F}^{n}(k) \right ]\\ 
\mathbb{E}\left[ \mathscr{F}(k) \right ] & \mathbb{E}\left[ \mathscr{F}^{2}(k) \right ] & \cdots & \mathbb{E}\left[ \mathscr{F}^{n+1}(k) \right ]\\ 
\vdots & \vdots & \ddots & \\ 
\mathbb{E}\left[ \mathscr{F}^{n}(k) \right ] & \mathbb{E}\left[ \mathscr{F}^{n+1}(k) \right ] &  & \mathbb{E}\left[ \mathscr{F}^{2n}(k) \right ]
\end{bmatrix} \succeq 0.
\end{aligned}
\label{HankelFk}
\end{equation}

Since each $\mathbb{E}\left[\mathscr{F}^i(k)\right]$ can be written as a function of $\mathbb{E}\left[ \tilde{u}^{\ell}(k) \right]$ and $\mathbb{E}\left[ w^{\ell}(k) \right]$ for $\ell \in \mathbb{N}_{1}^{2n}$ by \eqref{Expecttildeul}, \eqref{HankelFk} constitutes a constraint on $\mathbb{E}\left[ \tilde{u}^{\ell}(k) \right]$ for $\ell \in \mathbb{N}_{1}^{2n}$.

For the general Fokker-Planck control problem, we have the constraint on $\tilde{u}(k)$ proposed above, which leads to a constraint on the feasible desired terminal distribution. This constraint aligns intuitively with the nature of the Fokker-Planck equation, since it characterizes the diffusion process, which means that the uncertainty of the system state keeps increasing throughout the time evolution of the Fokker-Planck equation. Denote the $i_\text{th}$ order power moment of the desired terminal distribution as $\mathcal{M}^{i}_{\tau}$, namely
$$
\mathcal{M}^{i}_{\tau} = \int_{\mathbb{R}}x^{i}\tau(x)dx.
$$
As to verify the feasibility of a given moment sequence $\left\{ \mathcal{M}^{i}_{\tau}: i \in 
\mathbb{N}_{1}^{2n}\right\}$, we first minimize the uncertainty of the system state, namely to choose $c(k) = 1$. Then we have $x(k+1) = \tilde{u}(k)$, i.e.,
$$
\mathbb{E}\left[ \tilde{u}^{\ell}(k) \right] = \mathcal{M}^{\ell}_{\tau},
$$
for $\ell \in \mathbb{N}_{1}^{2n}$. Next, by substituting the corresponding $\left\{\mathbb{E}\left[ \tilde{u}^{\ell}(k) \right]: \ell \in \mathbb{N}_{1}^{2n}\right\}$ into \eqref{Expecttildeul}, we obtain the power moments $\left\{ \mathbb{E}\left[\mathscr{F}^{\ell}(k) \right]: \ell \in \mathbb{N}_{1}^{2n}\right\}$. Then if 
the Hankel matrix $H_{\mathscr{F}}(k)$ is positive semidefinite, the terminal distribution $\tau(x)$ is a reachable one. 

Moreover, we note that this result can not only be applied to the terminal desired distribution, but also for the distributions of all system states $x(k): k \in \mathbb{N}_{1}^{K}$. We conclude the following proposition.

\begin{proposition}
Given the Fokker-Planck equation \eqref{syseq} with control input \eqref{control}, the reachability condition of a distribution of system state $p_{k}(x)$ is that given
$$
\mathbb{E}\left[ \tilde{u}^{\ell}(k) \right] = \mathbb{E}\left[ x^{\ell}(k) \right] = \int_{\mathbb{R}}x^{\ell}p_{k}(x)dx,
$$
the Hankel matrix of $\mathscr{F}(k)$, of which the entries are calculated by \eqref{Expecttildeul}, is positive semidefinite. The corresponding state of the moment system $\mathscr{X}(k)$ is called a reachable state of the moment system.
\label{Proposition31}
\end{proposition}

\section{The control problem}
By Proposition \ref{Proposition31}, one would be able to ensure the reachability of a given system state, given the Fokker-Planck equation. We first consider determining the states of the moment system, namely $\left\{ \mathscr{X}(k): k \in \mathbb{N}_{0}^{K-1} \right\}$. For the general discrete-time Liouville control problem, suboptimal control schemes by convex optimization are proposed in \cite{wu2023general}. However, for the general Fokker-Planck problem, we are not always able to ensure that the state of the moment system $\mathscr{X}(k)$ satisfies the condition in Proposition \ref{Proposition31}, when the additive noise $w(k)$ is considerably large. In the literature, the trajectory of the system state is prespecified \cite{sivaramakrishnan2022distribution}. Considering this problem, we use convex optimization method to determine $\left\{ \mathscr{X}(k): k \in \mathbb{N}_{0}^{K-1} \right\}$ if the noise is not a dominating factor of the Fokker-Planck equation, and obtain $\left\{ \mathscr{X}(k): k \in \mathbb{N}_{0}^{K-1} \right\}$ heuristically if $w(k)$ is relatively large.

So then the challenge is now to determine the parameters $\left\{ c(k): k \in \mathbb{N}_{0}^{K-1} \right\}$ and the random variables $\left\{ \mathscr{F}(k): k \in \mathbb{N}_{0}^{K-1} \right\}$. To determine the latter one is equivalent to determine $\left\{ \tilde{u}(k): k \in 
\mathbb{N}_{0}^{K-1} \right\}$. Since we consider the moment system at this time being, the problem becomes determining the power moments up to order $2n$ of each $\tilde{u}(k)$ for $k \in \mathbb{N}_{0}^{K-1}$. Define the set of all reachable states $\mathscr{X}(k)$ of the moment system, given the additive noise $w(k)$, as $\Xi \left( w(k) \right)$. The optimization problem for the determination of $\left\{ c(k): k \in \mathbb{N}_{0}^{K-1} \right\}$ and $\left\{ \mathscr{F}(k): k \in \mathbb{N}_{0}^{K-1} \right\}$ can then be formulated as
\begin{equation}
\begin{aligned}
& \min_{c(k)}\mathbb{E}\left[ \left( -c(k)a(k)x(k) + \tilde{u}(k) \right)^{2} \right],\\
\mathrm{s.t.} \ & \tilde{\mathscr{U}}(k) = \mathscr{X}(k+1) - \tilde{\mathscr{A}}(\tilde{\mathscr{U}}(k))\mathscr{X}(k),\\
& H_{u}(k) \succeq 0, \\
& \mathscr{X}(k) \in \Xi \left( w(k) \right), \\
& 0 \leq c(k) \leq 1.
\end{aligned}
\label{optimization1}
\end{equation}
We first need to prove the existence and uniqueness of solution to the optimization problem.
\begin{theorem}
    The optimization problem \eqref{optimization1} is convex, and there exists a solution to this problem.
\label{Theorem41}
\end{theorem}
\begin{proof}
    To prove this theorem, we need to perform the following three proofs: (1) the cost function is convex; (2) the set of all feasible $c(k)$ is convex; (3) there exists a solution to \eqref{optimization1}.

    (1) The second order power moment of $u(k)$ reads
\begin{equation*}
\begin{aligned}
    & \mathbb{E}\left[ u^{2}(k) \right]\\
    = & c^{2}(k)\mathbb{E}\left[a^{2}(k)\right]\mathbb{E}\left[x^{2}(k)\right] + \mathbb{E}\left[\tilde{u}^{2}(k)\right]\\ 
    - & 2c(k)\mathbb{E}\left[a(k)\right]\mathbb{E}\left[x(k)\right]\mathbb{E}\left[\tilde{u}(k)\right].
\end{aligned}
\end{equation*}
The second order derivative is
$$
\frac{\mathrm{d}^{2}\mathbb{E}\left[ u^{2}(k) \right]}{\mathrm{d}c(k)^{2}} = \mathbb{E}\left[ a^{2}(k)
 \right]\mathbb{E}\left[ x^{2}(k) \right] \geq 0.
$$
Therefore, the cost function is convex.

(2) We need to prove that the set of all $c(k)$ satisfying all the four conditions in \eqref{optimization1} is a convex set.  
Because of the first and the third conditions in 
\eqref{optimization1}, there exists at least a $c(k)$ satisfying \eqref{momentsyseq}, and then \eqref{newsyseq}. Assume that $c_{1}(k), c_{2}(k)$ satisfy \eqref{newsyseq}, i.e.,
$$
x(k+1) = \left( 1 - b(k)c_{1}(k) \right)a(k)x(k) + \tilde{u}_{1}(k),
$$
and
$$
x(k+1) = \left( 1 - b(k)c_{2}(k) \right)a(k)x(k) + \tilde{u}_{2}(k).
$$
Let $0 \leq \lambda \leq 1$, we have that
$$
\begin{aligned}
& \lambda x(k+1) + (1 - \lambda) x(k+1)\\
= & \left( 1 - b(k)\left(\lambda c_{1}(k) + (1 - \lambda) c_{2}(k)\right) \right)a(k)x(k)\\
+ & \lambda \tilde{u}_{1}(k) + (1 - \lambda) \tilde{u}_{2}(k).
\end{aligned}
$$
Therefore, $\lambda c_{1}(k) + (1 - \lambda) c_{2}(k)$ is also a feasible solution. We can then conclude that the set of all feasible $c(k)$ is convex. However, we still need to prove the set is closed.

(3) Since the Hankel matrix $H_{u}(k)$ is a continuous matrix function of $c(k)$, and $H_{u}(k) \succeq 0$, the set of all feasible $c(k)$ is a compact convex set, which is a subset of $\left[ 0, 1 \right]$. It proves the existence of solution to the optimization problem \eqref{optimization1}.

\end{proof}

By Theorem \ref{Theorem41}, one is able to obtain $\left\{ c(k): k \in \mathbb{N}_{0}^{K-1} \right\}$, $\left\{ \tilde{\mathscr{U}}(k): k \in \mathbb{N}_{0}^{K-1} \right\}$, and then $\left\{\mathbb{E}\left[\mathscr{F}^i(k)\right]: i \in \mathbb{N}_{1}^{2n}\right\}$ for $k \in \mathbb{N}_{0}^{K-1}$. Now the problem comes to determining $\mathscr{F}(k)$ by the power moments for each $k \in \mathbb{N}_{0}^{K-1}$.

This is an ill-posed problem, which doesn't have a unique solution. We call it the realization of the controller. In this paper, we follow the treatment in our previous works \cite{wu2022group, wu2023general}. We ignore time step $k$ if there is no ambiguity.

Define the Kullback-Leibler divergence as
\begin{equation}
\mathbb{K} \mathbb{L}(r \| p)=\int_{\mathbb{R}} r(\mathscr{F}) \log \frac{r(\mathscr{F})}{p(\mathscr{F})} \mathrm{d}\mathscr{F}
\label{KLD}
\end{equation} 
where $r$ is a prescribed reference distribution in $\mathcal{P}_{2n}$. Define
$$
G(\mathscr{F})= \begin{bmatrix}
1 & \mathscr{F} & \cdots & \mathscr{F}^{n-1} & \mathscr{F}^{n},
\end{bmatrix}^{\intercal}
$$
and the linear integral operator $\Omega$ as
$$
\Omega: p(\mathscr{F}) \mapsto H_{\mathscr{F}}(k)=\int_{\mathbb{R}} G(\mathscr{F}) p(\mathscr{F}) G^{\intercal}(\mathscr{F}) \mathrm{d}\mathscr{F},
$$
where $p(\mathscr{F})$ is defined on the space $\mathcal{P}_{2n}$. It is obvious that $\mathcal{P}_{2n}$ is convex, and so is the range of $\Omega$, namely $\operatorname{range}(\Omega)=\Omega\mathcal{P}_{2n}$. Then the result is as follows.

Provided with a reference distribution $r \in \mathcal{P}_{2n}$ and a Hankel matrix $H_{\mathscr{F}}(k) \succ 0$, there is a unique $\hat{p} \in \mathcal{P}_{2n}$, which has the form
\begin{equation}
\hat{p}=\frac{r}{G^{\intercal} \hat{\Lambda} G}, 
\label{hatp}
\end{equation}
that minimizes \eqref{KLD} subject to $\Omega(\hat{p})=H_{\mathscr{F}}$. Here $\hat{\Lambda}$ is the unique solution to the problem of minimizing
\begin{equation}
\mathbb{J}_{r}(\Lambda):=\operatorname{tr}(\Lambda H_{\mathscr{F}})-\int_{\mathbb{R}} r(\mathscr{F}) \log \left[G(\mathscr{F})^{\intercal} \Lambda G(\mathscr{F})\right] \mathrm{d}\mathscr{F}
\label{Jr}
\end{equation}
over all $\Lambda \in \mathcal{L}_{+}$, which is defined as
$$
\mathcal{L}_{+}:=\left\{\Lambda \in \operatorname{range}(\Gamma) \mid G(\mathscr{F})^{\intercal} \Lambda G(\mathscr{F})>0, \mathscr{F} \in \mathbb{R}\right\}.
$$ 
$\operatorname{tr}(\cdot)$ denotes the trace of a matrix.

By this result, one would be able to turn the obtained power moments of $\mathscr{F}$ to the infinite-dimensional space $\mathcal{P}_{2n}$. And we conclude the results in the following Algorithm \ref{alg:1}.
\begin{algorithm}
    \caption{General discrete-time Fokker-Planck control by power moments.}
    \label{alg:1}
    \begin{algorithmic}[1]
        \Require The maximal time step $K$; the distribution of the system parameter $\left\{a(k), b(k): k\in \mathbb{N}_{0}^{K-1} \right\}$; a sample $x(0)$ drawn from the initial distribution $p_{0}(x)$; the specified terminal distribution $\tau(x)$.
        \Ensure The controls $u(k)$, $k = 0, \cdots, K-1$.
        \State $k \Leftarrow 0$
        \State Determine the states of the moment system $\left\{\mathscr{X}(k): k \in \mathbb{N}_{0}^{K}\right\}$ by convex optimization \cite{wu2023general} or by heuristic methods.
    \While{$0 \leq k < K$}
        \State Do optimization \eqref{optimization1}, obtain the $c(k)$ and the corresponding $\tilde{\mathscr{U}}(k)$ which minimize the convex cost function
        \State Calculate $\left\{\mathbb{E}\left[\mathscr{F}^i(k)\right]: i \in \mathbb{N}_{1}^{2n}\right\}$ by \eqref{Expecttildeul} using $\tilde{\mathscr{U}}(k)$.
        \State Optimize the cost function \eqref{Jr} and obtain the analytic estimate of the distribution $\hat{p}_{k}(\mathscr{F})$.
        \State Draw a sample, namely $\mathscr{F}(k)$, from $\hat{p}_{k}(\mathscr{F})$. 
        \State Calculate the control input $u(k)$ by \eqref{control} with the optimal $c(k)$ obtained in Step 4.
        \State $k \Leftarrow k+1$
    \EndWhile
    \end{algorithmic}
\end{algorithm}

\section{Numerical examples}
In this section, we will simulate two examples of the first-order discrete-time Fokker-Planck control problem. For each example, we perform $2000$ Monte-Carlo simulations. Then the histogram of the $2000$ terminal system states are compared to the desired distribution $\tau(x)$ to validate the proposed algorithm. 

In Example 1, we aim at steering a Gaussian distribution to a mixture of Gaussian distributions in $4$ steps. The system parameters $a(k) = 0.5$ and $b(k) = 0.8$ for $k \in \mathbb{N}_{0}^{3}$. The initial Gaussian distribution is chosen as 
\begin{equation}
     q_{0}(x) = \frac{1}{\sqrt{2\pi}}e^{\frac{x ^{2}}{2}},
\label{q01}
\end{equation}
and the terminal distribution is specified as
\begin{equation}
     \tau(x) = \frac{0.3}{\sqrt{2\pi} \cdot 2}e^{\frac{(x + 2) ^{2}}{2 \cdot 2^{2}}} + \frac{0.7}{\sqrt{2\pi} \cdot 2}e^{\frac{(x - 2) ^{2}}{2 \cdot 2^{2}}}. 
\label{qt1}
\end{equation}
The distribution of the white noise is chosen as $h_{k}(w) = \mathcal{N}(0, 1)$. The reference distribution $r(x)$ for the realization of $\mathscr{F}(k)$ is chosen as
$$
r(x) = \mathcal{N}(\mathbb{E}\left[ \mathscr{F}(k) \right], \mathbb{E}\left[ \mathscr{F}^{2}(k) \right]).
$$
In this example, we propose to consider the goal of optimization to be the maximal smoothness of state transition \cite{wu2023general}, i.e., 
$$
\min_{\mathscr{X}(1), \cdots, \mathscr{X}(K-1)} \sum_{k = 0}^{K-1}\left( \mathscr{X}(k+1) - \mathscr{X}(k) \right)^{T}\left( \mathscr{X}(k+1) - \mathscr{X}(k) \right),
$$
of which the solution is
\begin{equation}
\mathscr{X}(k) = \frac{K - k}{K} \mathscr{X}(0) + \frac{k}{K} \mathscr{X}(K).
\end{equation}
The system states of the moment system at time steps $\mathbb{N}_{1}^{K-1}$ are then determined.

\begin{figure}[htbp]
\centering
\includegraphics[scale=0.38]{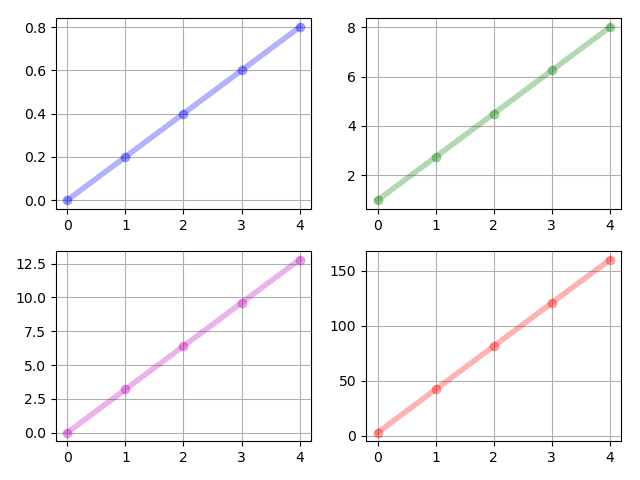}
\centering
\caption{$\mathscr{X}(k)$ at time steps $k = 0, 1, 2, 3, 4$. The upper left figure shows $\mathbb{E}\left[ x(k)\right]$. The upper right one shows $\mathbb{E}\left[ x^{2}(k)\right]$. The lower left one shows $\mathbb{E}\left[ x^{3}(k)\right]$ and the lower right one shows $\mathbb{E}\left[ x^{4}(k)\right]$.}
\label{fig1}
\end{figure}

\begin{figure}[htbp]
\centering
\includegraphics[scale=0.38]{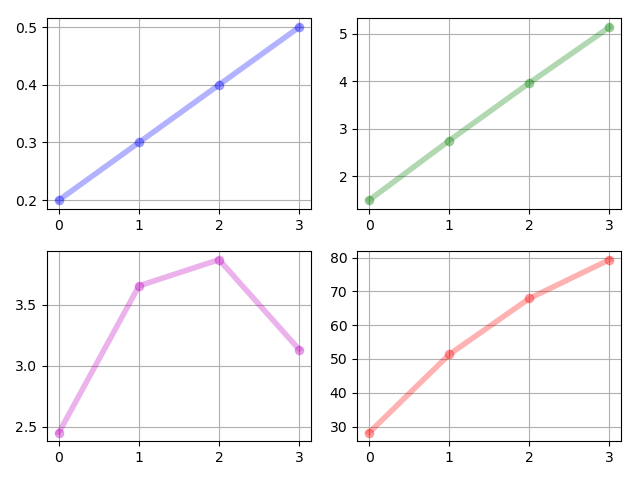}
\centering
\caption{The upper left figure shows $\mathbb{E}\left[ \mathscr{F}(k)\right]$. The upper right one shows $\mathbb{E}\left[ \mathscr{F}(k)^{2}\right]$. The lower left one shows $\mathbb{E}\left[\mathscr{F}(k)^{3}\right]$ and the lower right one shows $\mathbb{E}\left[ \mathscr{F}(k)^{4}\right]$.}
\label{fig2}
\end{figure}

\begin{figure}[htbp]
\centering
\includegraphics[scale=0.38]{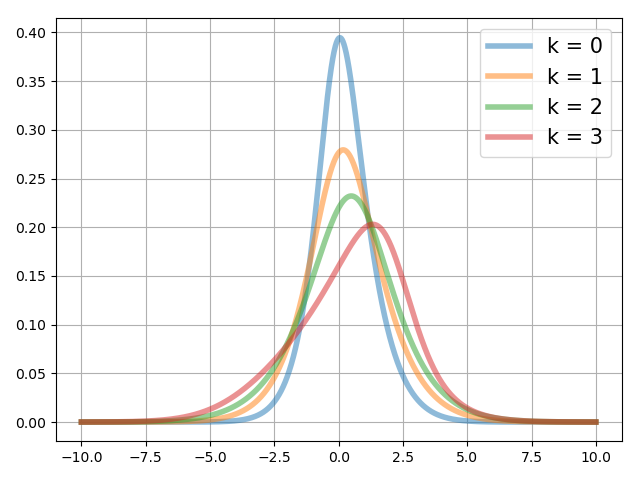}
\centering
\caption{Realized control inputs $\mathscr{F}(k)$ by $\left\{ \mathbb{E}\left[ \mathscr{F}^{i}(k) \right] : i \in \mathbb{N}_{1}^{2n}\right\}$ for $k = 0, 1, 2, 3$, which are obtained by our proposed control scheme.}
\label{fig3}
\end{figure}

\begin{figure}[htbp]
\centering
\includegraphics[scale=0.38]{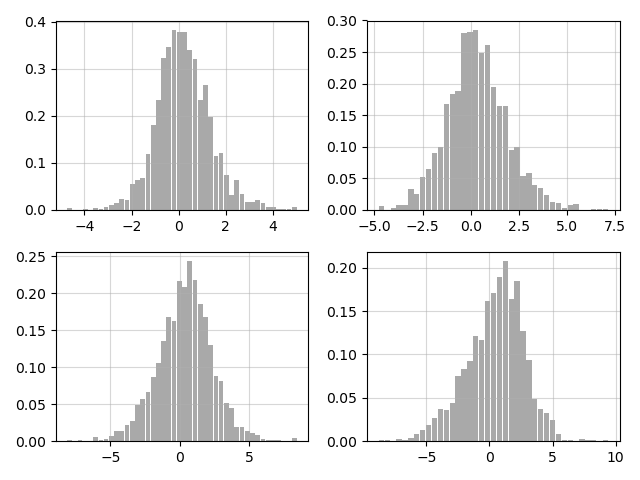}
\centering
\caption{The histograms of $u(k)$ at time step $k$ for the $2000$ Monte-Carlo simulations. The upper left and right figures are $u(0)$ and $u(1)$. The lower left and right figures are $u(2)$ and $u(3)$ respectively.}
\label{fig4}
\end{figure}

\begin{figure}[htbp]
\centering
\includegraphics[scale=0.38]{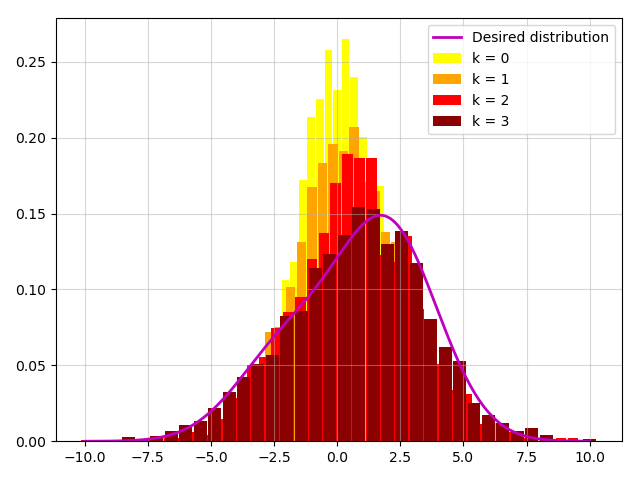}
\centering
\caption{The histograms of the system states $x(k)$ of $2000$ Monte-Carlo simulations at time step $k \in \mathbb{N}_{0}^{3}$, together with the desired terminal distribution $\tau(x)$. We note that the histogram of system state at step $k = 3$ is very close to $\tau(x)$.}
\label{fig5}
\end{figure}

The simulation results are given in Figure \ref{fig1} to Figure \ref{fig5}. The system states of the moment system, i.e., $\mathscr{X}(k)$ for $k \in \mathbb{N}_{0}^{3}$ are given in Figure \ref{fig1}. The power moments of $\mathscr{F}(k)$, i.e., $\left\{ \mathbb{E}\left[  \mathscr{F}^{i}(k) \right]: i \in \mathbb{N}_{1}^{2n} \right\}$ for $k \in \mathbb{N}_{0}^{3}$ are given in Figure \ref{fig2}. We note that by our proposed algorithm, $H_{x}(k)$ and $H_{\mathscr{F}}(k)$ are both positive definite, which makes it feasible for us to realize $\left\{\mathscr{F}(k): k \in \mathbb{N}_{0}^{3} \right\}$. The realized $\left\{\mathscr{F}(k): k \in \mathbb{N}_{0}^{3} \right\}$ are given in Figure \ref{fig3}. The optimal $c(k) = 0$ for $k = 0, \cdots, 3$. In Figure \ref{fig4}, the histograms of $\left\{u(k): k \in \mathbb{N}_{0}^{3}\right\}$ are provided. Histograms of system states $\left\{x(k): k \in \mathbb{N}_{0}^{3}\right\}$ are depicted in Figure \ref{fig5}. We observe that the transition of the histogram of the system state is very smooth. Moreover, the histogram of the terminal state of $2000$ Monte-Carlo simulations is very close to the desired terminal distribution $\tau(x)$. which validates the satisfactory performance of our proposed control scheme for the discrete-time first-order Fokker-Planck control problem together with the smooth state transition.

In Example 2, the task is to steer a Gaussian distribution to a mixture of generalized logistic distributions. The distribution of the initial system state is \eqref{q01}, and the desired terminal distribution is
\begin{equation}
     \tau(x) = \frac{0.4 \cdot 2e^{-x}}{(1 + e^{-x})^{3}} + \frac{0.6 \cdot 3e^{-x - 2}}{(1 + e^{-x - 2})^{4}}. 
\label{qt2}
\end{equation}

All other system settings are identical to that in Example 1. The states of the moment system, namely $\left\{ \mathscr{X}(k): k \in \mathbb{N}_{0}^{3}\right\}$, are obtained by the maximal smoothness of the state transition, following the treatment in Example 1. 

The simulation results have been visually depicted in Figures \ref{fig1} through \ref{fig5}. Notably, the optimal value for \(c(k)\) is \(0\) for \(k = 0, \cdots, 3\). Our observations highlight the remarkably smooth transition evident in the histograms of the system states. Furthermore, it is noteworthy that the histogram of the terminal state from \(2000\) Monte-Carlo simulations closely approximates the desired terminal distribution \(\tau(x)\). This compelling example serves to affirm the efficacy of the proposed control scheme in tackling discrete-time first-order Fokker-Planck control problems with non-Gaussian distributions.

\begin{figure}[htbp]
\centering
\includegraphics[scale=0.38]{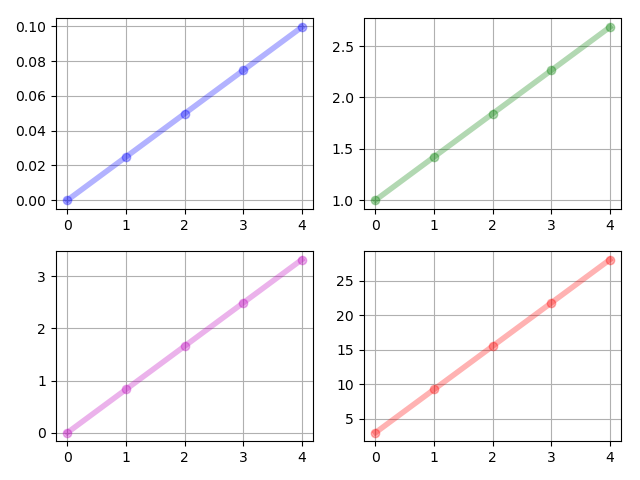}
\centering
\caption{$\mathscr{X}(k)$ at time steps $k = 0, 1, 2, 3, 4$. The upper left figure shows $\mathbb{E}\left[ x(k)\right]$. The upper right one shows $\mathbb{E}\left[ x^{2}(k)\right]$. The lower left one shows $\mathbb{E}\left[ x^{3}(k)\right]$ and the lower right one shows $\mathbb{E}\left[ x^{4}(k)\right]$.}
\label{fig6}
\end{figure}

\begin{figure}[htbp]
\centering
\includegraphics[scale=0.38]{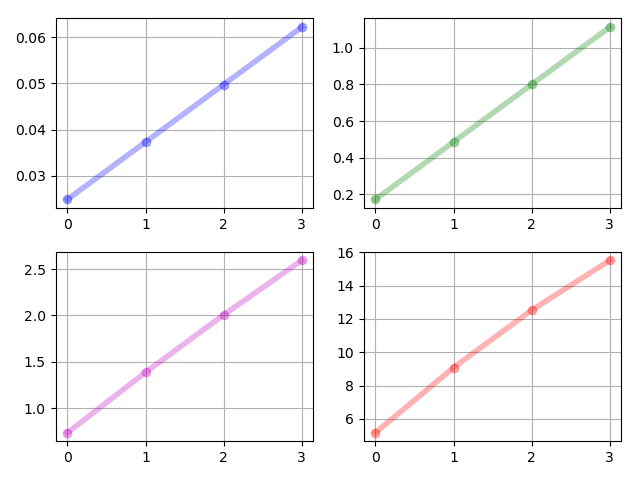}
\centering
\caption{The upper left figure shows $\mathbb{E}\left[ \mathscr{F}(k)\right]$. The upper right one shows $\mathbb{E}\left[ \mathscr{F}(k)^{2}\right]$. The lower left one shows $\mathbb{E}\left[\mathscr{F}(k)^{3}\right]$ and the lower right one shows $\mathbb{E}\left[ \mathscr{F}(k)^{4}\right]$.}
\label{fig7}
\end{figure}

\begin{figure}[htbp]
\centering
\includegraphics[scale=0.38]{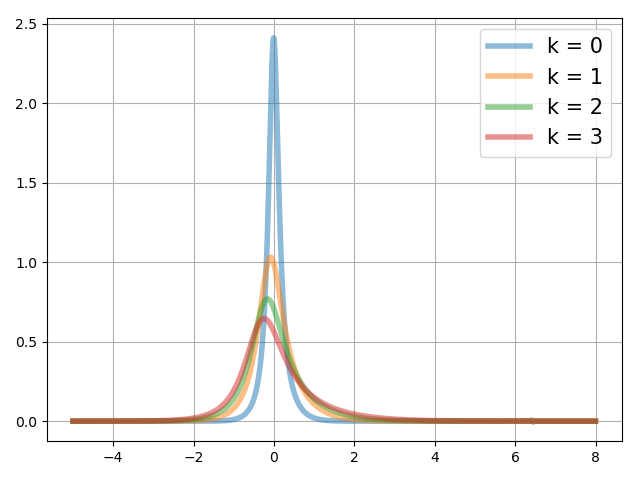}
\centering
\caption{Realized control inputs $\mathscr{F}(k)$ by $\left\{ \mathbb{E}\left[ \mathscr{F}^{i}(k) \right] : i \in \mathbb{N}_{1}^{2n}\right\}$ for $k = 0, 1, 2, 3$, which are obtained by our proposed control scheme.}
\label{fig8}
\end{figure}

\begin{figure}[htbp]
\centering
\includegraphics[scale=0.38]{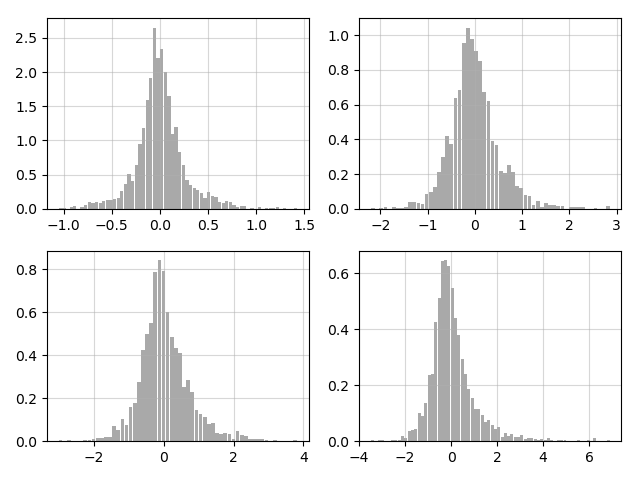}
\centering
\caption{The histograms of $u(k)$ at time step $k$ for the $2000$ Monte-Carlo simulations. The upper left and right figures are $u(0)$ and $u(1)$. The lower left and right figures are $u(2)$ and $u(3)$ respectively.}
\label{fig9}
\end{figure}

\begin{figure}[htbp]
\centering
\includegraphics[scale=0.38]{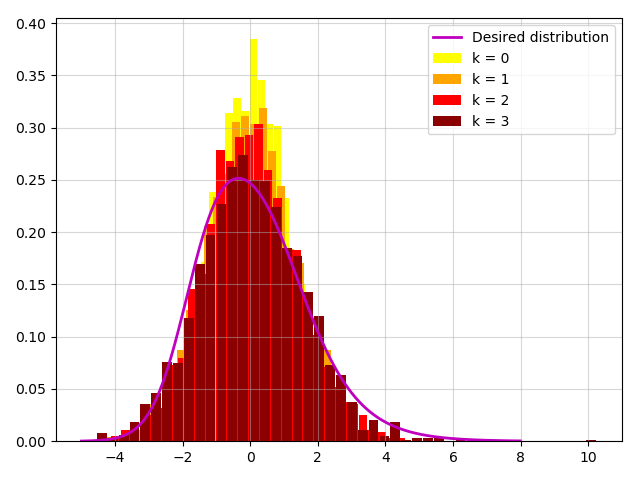}
\centering
\caption{The histogram of the system states $x(k)$ of $2000$ Monte-Carlo simulations at time step $k \in \mathbb{N}_{0}^{3}$, together with the desired terminal distribution $\tau(x)$. The histogram of system state at step $k = 3$ is very close to $\tau(x)$.}
\label{fig10}
\end{figure}

\section{A concluding remark}
We introduce a systematic methodology for addressing the general discrete-time Fokker-Planck control problem, extending beyond the realm of Gaussian distributions for initial and terminal system state distributions. Instead of imposing Gaussian constraints, we require the initial and terminal distributions to possess the first $2n$ orders of power moments, where $n$ is a flexible positive integer. Through this approach, we transform the original problem, which resides inherently in infinite dimensions, into a finite-dimensional context using the power moments. Our proposed method diverges from conventional deterministic feedback control methods, which are proven in \cite{elamvazhuthi2018optimal} to be infeasible for this intricate problem.

Instead, we devise an innovative control strategy amalgamating a feedback component and a random variable that serves as a transition kernel for the underlying Markovian process. We establish conditions under which this control approach is viable, contingent upon the discrete-time Fokker-Planck equation and the accompanying additive system noise $w(k)$. To optimize system parameters, we present a convex optimization scheme, together with the proof to the existence of solution.

A realization of the transition kernel from its power moments to the infinite-dimensional space $\mathcal{P}_{2n}$ is proposed. Algorithm 1 furnishes a comprehensive procedure for tackling the general discrete-time Fokker-Planck control problem. Notably, we validate our proposed control strategy through simulation of two numerical examples, of which the terminal distributions are a mixture of two Gaussian distributions and a mixture of two non-Gaussian distributions, respectively. This substantiates the effectiveness of our approach.

In forthcoming research, we aspire to broaden the scope of this paper's findings to encompass more intricate multivariate discrete-time Fokker-Planck control challenges. This extension aims to enhance the applicability of our algorithm in real-world engineering scenarios, although it constitutes a non-trivial task. The parameterization of the multivariate transition kernel's realization presents not only a challenge for controller design but also a significant inquiry within the realm of the Hamburger (Hausdorff) moment problem.

\bibliographystyle{plain}
\bibliography{autosam}

\begin{thebibliography}{10}

\bibitem{bakolas2019dynamic}
Efstathios Bakolas.
\newblock Dynamic output feedback control of the liouville equation for
  discrete-time siso linear systems.
\newblock {\em IEEE Transactions on Automatic Control}, 64(10):4268--4275,
  2019.

\bibitem{balci2020covariance}
Isin~M Balci and Efstathios Bakolas.
\newblock {Covariance steering of discrete-time stochastic linear systems based
  on Wasserstein distance terminal cost}.
\newblock {\em IEEE Control Systems Letters}, 5(6):2000--2005, 2020.

\bibitem{balci2021covariance}
Isin~M Balci and Efstathios Bakolas.
\newblock Covariance control of discrete-time gaussian linear systems using
  affine disturbance feedback control policies.
\newblock In {\em 2021 60th IEEE Conference on Decision and Control (CDC)},
  pages 2324--2329. IEEE, 2021.

\bibitem{balci2023covariance}
Isin~M Balci and Efstathios Bakolas.
\newblock Covariance steering of discrete-time linear systems with mixed
  multiplicative and additive noise.
\newblock In {\em 2023 American Control Conference (ACC)}, pages 2586--2591.
  IEEE, 2023.

\bibitem{biswal2021decentralized}
Shiba Biswal, Karthik Elamvazhuthi, and Spring Berman.
\newblock Decentralized control of multiagent systems using local density
  feedback.
\newblock {\em IEEE Transactions on Automatic Control}, 67(8):3920--3932, 2021.

\bibitem{blaquiere1992controllability}
A~Blaquiere.
\newblock {Controllability of a Fokker-Planck equation, the Schr{\"o}dinger
  system, and a related stochastic optimal control (revised version)}.
\newblock {\em Dynamics and Control}, 2(3):235--253, 1992.

\bibitem{brockett2012control}
R~Brockett, J~Le~Rousseau, O~Glass, and E~Zuazua.
\newblock Control of partial differential equations.
\newblock In {\em New York}. Springer, 2012, 2012.

\bibitem{chen2015optimal}
Yongxin Chen, Tryphon~T Georgiou, and Michele Pavon.
\newblock {Optimal steering of a linear stochastic system to a final
  probability distribution, Part I}.
\newblock {\em IEEE Transactions on Automatic Control}, 61(5):1158--1169, 2015.

\bibitem{chen2015optimal2}
Yongxin Chen, Tryphon~T Georgiou, and Michele Pavon.
\newblock {Optimal steering of a linear stochastic system to a final
  probability distribution, Part II}.
\newblock {\em IEEE Transactions on Automatic Control}, 61(5):1170--1180, 2015.

\bibitem{chen2016optimal}
Yongxin Chen, Tryphon~T Georgiou, and Michele Pavon.
\newblock Optimal transport over a linear dynamical system.
\newblock {\em IEEE Transactions on Automatic Control}, 62(5):2137--2152, 2016.

\bibitem{chen2018optimal}
Yongxin Chen, Tryphon~T Georgiou, and Michele Pavon.
\newblock {Optimal steering of a linear stochastic system to a final
  probability distribution—part III}.
\newblock {\em IEEE Transactions on Automatic Control}, 63(9):3112--3118, 2018.

\bibitem{dirr2016controlling}
G~Dirr, U~Helmke, and M~Sch{\"o}nlein.
\newblock Controlling mean and variance in ensembles of linear systems.
\newblock {\em IFAC-PapersOnLine}, 49(18):1018--1023, 2016.

\bibitem{elamvazhuthi2018optimal}
Karthik Elamvazhuthi, Piyush Grover, and Spring Berman.
\newblock Optimal transport over deterministic discrete-time nonlinear systems
  using stochastic feedback laws.
\newblock {\em IEEE control systems letters}, 3(1):168--173, 2018.

\bibitem{filliger2008connection}
Roger Filliger, M-O Hongler, and Ludwig Streit.
\newblock Connection between an exactly solvable stochastic optimal control
  problem and a nonlinear reaction-diffusion equation.
\newblock {\em Journal of Optimization Theory and Applications}, 137:497--505,
  2008.

\bibitem{fuhrmann2015mathematics}
Paul~Abraham Fuhrmann, Uwe Helmke, et~al.
\newblock {\em The mathematics of networks of linear systems}, volume 150.
\newblock Springer, 2015.

\bibitem{goldshtein2017finite}
Maxim Goldshtein and Panagiotis Tsiotras.
\newblock Finite-horizon covariance control of linear time-varying systems.
\newblock In {\em 2017 IEEE 56th Annual Conference on Decision and Control
  (CDC)}, pages 3606--3611. IEEE, 2017.

\bibitem{li2010ensemble}
Jr-Shin Li.
\newblock Ensemble control of finite-dimensional time-varying linear systems.
\newblock {\em IEEE Transactions on Automatic Control}, 56(2):345--357, 2010.

\bibitem{li2007ensemble}
Jr-Shin Li and Navin Khaneja.
\newblock Ensemble control of linear systems.
\newblock In {\em 2007 46th IEEE conference on decision and control}, pages
  3768--3773. IEEE, 2007.

\bibitem{liu2022optimal}
Fengjiao Liu and Panagiotis Tsiotras.
\newblock Optimal covariance steering for continuous-time linear stochastic
  systems with additive generic noise.
\newblock {\em arXiv preprint arXiv:2206.11201}, 2022.

\bibitem{okamoto2018optimal}
Kazuhide Okamoto, Maxim Goldshtein, and Panagiotis Tsiotras.
\newblock Optimal covariance control for stochastic systems under chance
  constraints.
\newblock {\em IEEE Control Systems Letters}, 2(2):266--271, 2018.

\bibitem{pavon1991free}
Michele Pavon and Anton Wakolbinger.
\newblock On free energy, stochastic control, and schr{\"o}dinger processes.
\newblock In {\em Modeling, Estimation and Control of Systems with Uncertainty:
  Proceedings of a Conference held in Sopron, Hungary, September 1990}, pages
  334--348. Springer, 1991.

\bibitem{pilipovsky2023data}
Joshua Pilipovsky and Panagiotis Tsiotras.
\newblock Data-driven covariance steering control design.
\newblock {\em arXiv preprint arXiv:2303.17675}, 2023.

\bibitem{pra1990markov}
Paolo~Dai Pra and Michele Pavon.
\newblock {On the Markov processes of Schr{\"o}dinger, the Feynman-Kac formula
  and stochastic control}.
\newblock In {\em Realization and Modelling in System Theory: Proceedings of
  the International Symposium MTNS-89, Volume I}, pages 497--504. Springer,
  1990.

\bibitem{qi2013optimal}
Ji~Qi, Anatoly Zlotnik, and Jr-Shin Li.
\newblock Optimal ensemble control of stochastic time-varying linear systems.
\newblock {\em Systems \& Control Letters}, 62(11):1057--1064, 2013.

\bibitem{sivaramakrishnan2022distribution}
Vignesh Sivaramakrishnan, Joshua Pilipovsky, Meeko Oishi, and Panagiotis
  Tsiotras.
\newblock Distribution steering for discrete-time linear systems with general
  disturbances using characteristic functions.
\newblock In {\em 2022 American Control Conference (ACC)}, pages 4183--4190.
  IEEE, 2022.

\bibitem{wu2022density}
Guangyu Wu and Anders Lindquist.
\newblock Density steering by power moments.
\newblock {\em arXiv preprint arXiv:2211.02322}, 2022.

\bibitem{wu2022group}
Guangyu Wu and Anders Lindquist.
\newblock Group steering: Approaches based on power moments.
\newblock {\em arXiv preprint arXiv:2211.13370}, 2022.

\bibitem{wu2023general}
Guangyu Wu and Anders Lindquist.
\newblock General distribution steering: A sub-optimal solution by convex
  optimization.
\newblock {\em arXiv preprint arXiv:2301.06227}, 2023.

\end{thebibliography}

\begin{wrapfigure}{l}{20mm} 
\includegraphics[width=1in,height=1.25in,clip,keepaspectratio]{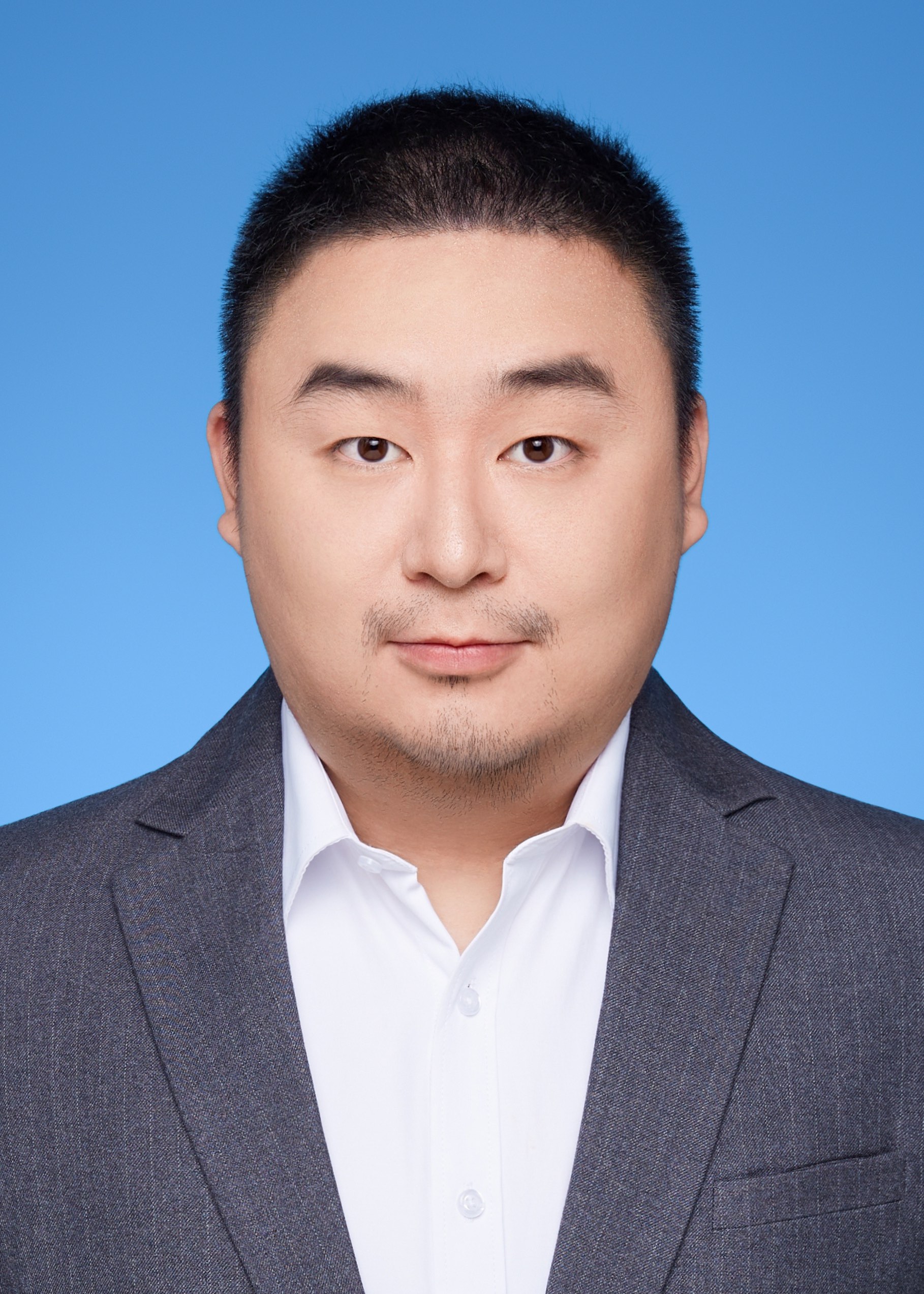}
  \end{wrapfigure}\par
  \textbf{Guangyu Wu} received the B.E. degree from Northwestern Polytechnical University, Xi’an, China, in 2013, and two M.S. degrees, one in control science and engineering from Shanghai Jiao Tong University, Shanghai, China, in 2016, and the other in electrical engineering from the University of Notre Dame, South Bend, USA, in 2018. 

He is currently pursuing the Ph.D. degree at Shanghai Jiao Tong University. His research interests are the moment problem and its applications to stochastic filtering, stochastic control, system identification and statistics. He is an active reviewer of journals and conferences including Automatica and International Conference on Machine Learning (ICML).

\begin{wrapfigure}{l}{25mm} 
\includegraphics[width=1in,height=1.25in,clip,keepaspectratio]{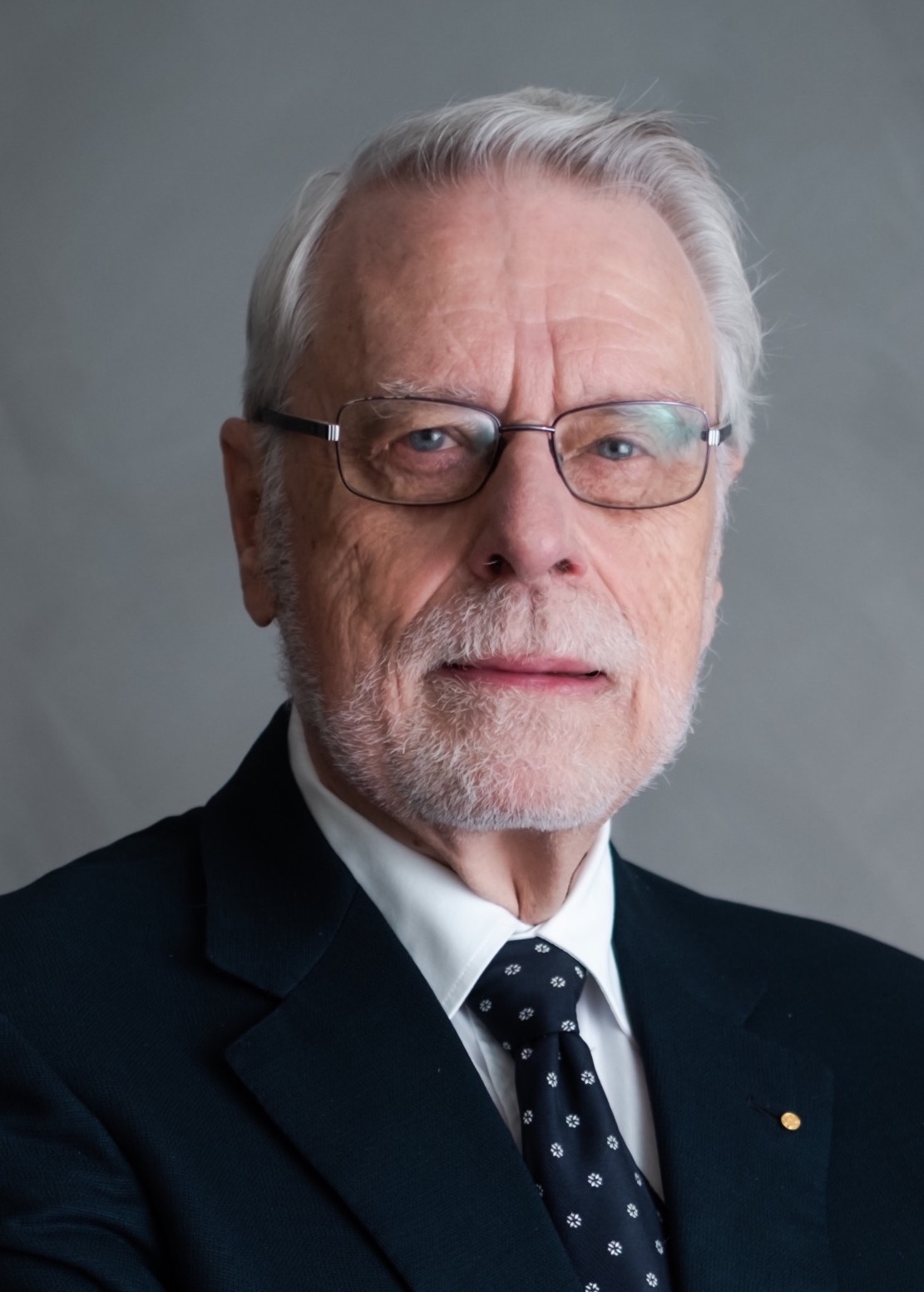}
  \end{wrapfigure}\par
  \textbf{Anders Lindquist} received the Ph.D. degree in optimization and systems theory from the Royal Institute of Technology, Stockholm, Sweden, in 1972, and an honorary doctorate (Doctor Scientiarum Honoris Causa) from Technion (Israel Institute of Technology) in 2010. 

He is currently a Zhiyuan Chair Professor at Shanghai Jiao Tong University, China, and Professor Emeritus at the Royal Institute of Technology (KTH), Stockholm, Sweden. Before that he had a full academic career in the United States, after which he was appointed to the Chair of Optimization and Systems at KTH.
Dr. Lindquist is a Member of the Royal Swedish Academy of Engineering Sciences, a Foreign Member of the Chinese Academy of Sciences, a Foreign Member of the Russian Academy of Natural Sciences, a Member of Academia Europaea (Academy of Europe), an Honorary Member the Hungarian Operations Research Society, a Fellow of SIAM, and a Fellow of IFAC. He received the 2003 George S. Axelby Outstanding Paper Award, the 2009 Reid Prize in Mathematics from SIAM, and the 2020 IEEE Control Systems Award, the
IEEE field award in Systems and Control.

\end{document}